\documentclass[letterpaper, 10 pt, conference]{ieeeconf}  

\IEEEoverridecommandlockouts                              
\title{\LARGE \bf
Nonlinear Fisher Particle Output Feedback Control and its application to Terrain Aided Navigation 
}

\usepackage[autostyle]{csquotes}
\usepackage{amsmath,verbatim,amssymb,amsfonts,amscd, graphicx}
\usepackage{mathabx}
\usepackage{breqn}
\usepackage{float}
\usepackage{graphics}
\usepackage{stmaryrd} 
\usepackage{textcomp}
\usepackage{setspace}
\usepackage{algorithm}
\usepackage[noend]{algpseudocode}

\author{Emilien Flayac, Karim Dahia, Bruno H\'eriss\'e, and Fr\'ed\'eric Jean
\thanks{}
\thanks{
        {\tt\small }}%
\thanks{
        {\tt\small}}%
}

\begin{document}

\maketitle
\thispagestyle{empty}
\pagestyle{empty}

\begin{abstract}
This paper presents state estimation and stochastic optimal control gathered in one global optimization problem generating  dual effect i.e. the control can improve the future estimation. As the optimal policy is impossible to compute, a sub-optimal policy that preserves this coupling is constructed thanks to the Fisher Information Matrix (FIM) and a Particle Filter. This method has been applied to the localization and guidance of a drone over a known terrain with height measurements only. The results show that the new method  improves the estimation accuracy compared to nominal trajectories. 
\end{abstract}


\section*{Introduction}

Stochastic optimal control problems with imperfect state information arise when an optimal control problem contains uncertainties on the dynamics and when its state is partially observed. These problems have many applications in chemistry \cite{bayard_implicit_2008}\cite{subramanian_economic_2015} and in the automotive industry \cite{blackmore_probabilistic_2010} for unmanned vehicles, for example.
The Dynamic Programming principle \cite{bertsekas_dynamic_2005} theoretically allows one to find the optimal controls looked for as policies due to the randomness of the problem. In addition, in such problems, as one has an access to the state of the system only through some observations, a state estimator is also needed as an as a function of them. 
In some problems, the observations depend on the control, it is then said that the control has a \textit{dual effect} \cite{feldbaum_optimal_1965}. It has a double role: it guides the system in a standard way and, at the same time,  it can also look for more information about the system because it influences the observations \cite{bar-shalom_dual_1974}.  

Optimal policies are often impossible to compute directly because of the curse of dimensionality. Thus many sub-optimal policies have been developed to approximate the optimal one. A sub-optimal policy can be designed to keep the property of dual effect. It is mostly done when the control problem is mixed with a \textit{parameter} estimation problem \cite{feldbaum_optimal_1965}. Indeed, these methods are applied when learning about an unknown parameter of a system helps guiding it. We present a problem where the dual effect is used to improve \textit{state} estimation. Particle approximations are then very promising techniques. Indeed, they are very efficient to approximate stochastic optimization problems or to  estimate the state of a system, even in presence of high uncertainties, high non-linearities and probability constraints. 

Particle approximations are widely used in robust control. In \cite{blackmore_chance-constrained_2011}, the planned trajectories consider uncertainties, obstacles or other probability constraints. Nevertheless, these methods do not include state estimation and do not compute control policies but control values. In \cite{copp_nonlinear_2014}, an optimization problem coupling state estimation by a \textit{Moving Horizon Estimation} (MHE) and control by \textit{Model Predictive Control} (MPC) is discussed but this problem does not include dual effect. In \cite{sehr_particle_2016} and \cite{stahl_pf-mpc:_2011}, a Particle Output MPC policy with a particle filter used for the estimation and inside the optimization problem is presented but, again, there is no coupling between the control and the future estimation. In \cite{hanssen_scenario_2015}, a dual controller based on a tree representation by particles is proposed. However, in the latter article, the particles inside the optimization problems are introduced by an Ensemble Kalman filter rather than with a Particle Filter. In \cite{bayard_implicit_2008}, an implicit dual controller is computed thanks to a particle-based policy iteration approximation but it is extremely costly in practice and is limited to finite control spaces. In \cite{subramanian_economic_2015}, an Output feedback method based on Unscented Kalman filter with a tree representation and measurements anticipation is proposed but the conditional probability density of the state is supposed to be gaussian at each time.

In this paper, we propose a particular stochastic optimization problem that merges state estimation and control. This problem makes explicitly appear dual effect additively in its cost which creates a coupling between the controls and the state estimators.

We also propose a sub-optimal policy of our new optimization problem based on two successive approximations. The first one consists in replacing a term by an equivalent and simpler one which maintains the coupling created by the dual effect. The second one is a particle approximation used, both inside the optimization problem to find the control, and outside it to estimate the state.

This paper is organized the following way: in section \ref{pb_opti}, we describe our new stochastic optimization problem and give a comparison with classical problems. In section \ref{tract_app},  we describe the approximation of our problem and compare it to existing ones. We also give an application of our method with numerical results.

\section{Setup of stochastic optimal control}\label{pb_opti}
\subsection{Optimization  problem coupling control and estimation} \label{sec_pce}
\subsubsection{Stochastic dynamics and observation equation }
We consider a discrete-time stochastic dynamical system whose state is a stochastic process ${(X_k)}_{k\in \llbracket 0,T \rrbracket}$ valued in $\mathbb{R}^n$ with $T\in \mathbb{N^*}$ which verifies $\forall k \in \llbracket 0,T-1 \rrbracket$:

\begin{align}
X_{k+1}&=f_k(X_{k},U_{k},{\xi}_{k}) \label{sys_dyn},\\
X_0&\sim p_0,\notag
\end{align}
where:
\begin{itemize}
{\setlength{\baselineskip}{1.2\baselineskip}
\item $p_0$ is a  probability density and $X_0\sim p_0$ means that $p_0$ is the probability law of $X_0$.
\item 
${(U_k)}_{k\in \llbracket 0,T \rrbracket}$ is a stochastic process such that $\forall k\in \llbracket 0,T-1\rrbracket$, $U_k$ is valued in $\mathbb{U}_k \subset \mathbb{R}^m$.
\item ${(\xi_k)}_{k\in \llbracket 0,T \rrbracket}$ is a stochastic process valued in  $\mathbb{R}^d$ which corresponds to the disturbances on the dynamics. We suppose that $\forall k\in \llbracket 0,T-1\rrbracket$, $\xi_k \sim p_{\xi_k}$, and that $\xi_k$ is independent of ${\xi_l}$ for $k\neq l$ and of $X_0$.
\par}
\item $\forall k \in \llbracket 0,T-1 \rrbracket$, $f_k : \mathbb{R}^n\times\mathbb{R}^m\times\mathbb{R}^d \longrightarrow \mathbb{R}^n$.
\end{itemize}
We also assume that the state of the system is available through some observations represented by  a stochastic process ${(Z_k)}_{k\in\llbracket 0,T \rrbracket}$ valued in $\mathbb{R}^d$  which verifies, $\forall k \in \llbracket 0,T \rrbracket$:
\begin{align}
Z_k=h_k(X_k,\eta_k),\label{eq_obs}
\end{align}
where:
\begin{itemize}
{\setlength{\baselineskip}{1.2\baselineskip}
\item ${(\eta_k)}_{k\in \llbracket 0,T \rrbracket}$ is a stochastic process valued in  $\mathbb{R}^q$ which corresponds to the disturbances on the observations. We suppose that $\forall k\in \llbracket 0,T\rrbracket$, $\eta_k \sim p_{\eta_k}$ and that $\eta_k$ is independent of $\xi_k$, $X_0$ and  ${\eta_l}$ for $k\neq l$.
\par}
\item $\forall k \in \llbracket 0,T \rrbracket$, $h_k : \mathbb{R}^n\times\mathbb{R}^q \longrightarrow \mathbb{R}^p$.
\end{itemize}
For $k\in \llbracket 0,T\rrbracket$, we define the \textit{information vector}  $I_k$ such as:
\vspace{-0.2cm}
\begin{align}
&I_0=Z_0,&I_{k+1}=(I_k,U_{k},Z_{k+1}).\label{inf_vec}
\end{align}
\subsubsection{Presentation of our new optimization problem}
As explained in \cite{bertsekas_dynamic_2005} and \cite{bayard_implicit_2008}, in stochastic control, one does not seek control \textit{values} like in deterministic control but \textit{policies} i.e. functions of a certain random variable. As $I_k$ gathers all the data available for the controller, $U_k$ will be looked for as a function of $I_k$.  Moreover, for the same reason about $I_k$, any estimator of $X_k$, denoted by $\widehat{X}_k$, will also be looked for as a function of $I_k$.
 Starting from this remark, $\forall k \in \llbracket 0,T-1 \rrbracket$, we define a \textit{generalized} control $V_k=(U_k,\widehat{X}_k)$ and $V_T=\widehat{X}_T$ that must verify:
 \begin{align}
   V_k&=(U_k,\widehat{X}_{k})=(\mu_k(I_k),\pi_k(I_k)), \label{func_constraint}\\
  V_T&=\widehat{X}_{T}=\pi_T(I_T),\notag
 \end{align}
 where $\mu_k$ maps an information vector $I_k$ to a control $U_k$ in the control space $\mathbb{U}_k$ and $\pi_k$ maps an information vector $I_k$ to an estimator $\widehat{X}_{k}$ in $\mathbb{R}^n$. Thus, minimizing over $(V_0,\dots,V_T)$ with the constraints (\ref{func_constraint}) is equivalent to directly minimizing over  $(\mu_0,\dots,\mu_{T-1})$ and $(\pi_0,\dots,\pi_{T})$.
 Finally, similarly to what is done in \cite{copp_nonlinear_2014}, we propose a stochastic optimization problem over the generalized control $V_k$ that mixes control and state estimation. In addition, in our proposed approach, $(U_0,\dots,U_{T-1})$ and $(\widehat{X}_0,\dots,\widehat{X}_T)$ are coupled, which means that the control $U_k$ can influence the future estimators $(\widehat{X}_{k+1},\dots,\widehat{X}_{T})$. In order to do this, we define  generalized integral costs ${(\tilde{g}_k)}_{k \in \llbracket 0,T-1 \rrbracket}$ and a generalized  final cost
$\tilde{g}_N$ such that, $\forall k \in \llbracket 0,T-1 \rrbracket$, each one can be decomposed in two terms as follow:
\begin{align}
 \tilde{g}_k(X_k,V_k,\xi_k) &={g}_k(X_k,U_k,\xi_k)+f_C\left(C_k\right)\label{gen_int_cost},\\
\tilde{g}_T(X_T,V_T) &={g}_T(X_T)+f_C\left(C_T\right) \label{gen_fin_cost},
\end{align}
with
\begin{align}
 C_k&=\mathbb{E}\left[ (X_k-\widehat{X}_k){(X_k-\widehat{X}_k)}^T\right],\label{cout_cov_quad}
\end{align}
where:
\begin{itemize}
\item $\forall k \in \llbracket 0,T-1 \rrbracket$, $g_k$: $\mathbb{R}^n\times\mathbb{R}^m\times\mathbb{R}^d \longrightarrow \mathbb{R}$ is a standard instantaneous cost and $g_N$:$\mathbb{R}^n \longrightarrow \mathbb{R}$ is a standard final cost. Here, standard means that these costs are a criterion of the system performance we want to optimize in the first place like a price or a distance for example.
\item $f_C$ : $S_n^{++}(\mathbb{R})  \longrightarrow \mathbb{R}$ is a cost on the covariance matrix of the estimator, defined in (\ref{cout_cov_quad}). $f_C$ can be seen as a measure of the estimation error.
\end{itemize}
Therefore, minimizing the costs defined in (\ref{gen_int_cost}) and (\ref{gen_fin_cost}) over  $(V_0,\dots,V_T)$ is  equivalent to looking for a compromise between control and state estimation. With (\ref{sys_dyn})-(\ref{cout_cov_quad}), we can define our generalized stochastic optimal control problem $(P_{CE})$ by: 
\begin{equation*}
\begin{array}{rrclcl}
\underset{\underset{\mu_{0},\dots, \mu_{T-1} }{\pi_{0},\dots, \pi_{T}}}{\text{min}} & \multicolumn{3}{l}{\mathbb{E}\left[\sum_{k=0}^{T-1} \tilde{g}_k(X_k,V_k,\xi_k) + \tilde{g}_T(X_T,V_T)\right]} \\
\textrm{s.t.} &\forall k  & \in & \llbracket 0,T-1 \rrbracket, \\
& X_{k+1} & = & f_k(X_{k},U_{k},{\xi}_{k}), \\
& Z_k & = & h_k(X_k,\eta_k), \\
& V_k & = &(\mu_k(I_k),\pi_k(I_k)), \\
&  Z_T & = & h(X_T,\eta_T), \\
& V_T & = & \pi_T(I_T). \\
\end{array}
\end{equation*}
With an appropriate choice of $f_C$, the terms $f_C\left(C_k\right)$ can force a coupling between $U_{k-1}$ and $\widehat{X}_k$ and in particular the control $U_{k-1}$ can force the state $X_k$ to reduce the error made by the estimator $\widehat{X}_k$. Eventually, the sum of those terms creates
a coupling between $U_{k-1}$ and $(\widehat{X}_k,\dots,\widehat{X}_T)$.
Still, $(P_{CE})$ is computationally intractable because  $(\pi_{0},\dots,\pi_T)$  and $(\mu_{0},\dots, \mu_{T-1})$ are extremely hard to compute due to the curse of dimensionality.
Moreover, if $f_C$ is not linear, classical Dynamic Programming cannot be applied. 
In the following, we show, as in \cite{copp_nonlinear_2014}, that $(P_{CE})$ is a combination of two types of problems: a classical stochastic optimal control problem without state estimation and a sequence of  state estimation problems with a a-priori-fixed control. 
\subsection{Link with classical stochastic optimal control }
If one chooses $f_C$ to be constant then only remains the minimization over $(\mu_0,\dots,\mu_{T-1})$ and one recovers a stochastic optimal control problem with imperfect state information, denoted by $(P_C)$:
\begin{equation*}
\begin{array}{rrclcc}
\displaystyle (P_C):\underset{\mu_{0},\dots, \mu_{T-1}}{\text{min}} & \multicolumn{3}{l}{\mathbb{E}\left[\sum_{k=0}^{T-1} g_k(X_k,U_k,\xi_k) + g_T(X_T)\right]} \\
\textrm{s.t.} & X_{k+1} & = & f_k(X_{k},U_{k},{\xi}_{k}), \\
& Z_k & = & h_k(X_k,\eta_k), \\
& U_k & = & \mu_k(I_k), \; \forall k \in \llbracket 0,T-1 \rrbracket. \\
\end{array}
\end{equation*}
As shown in \cite{bertsekas_dynamic_2005}, the optimal policies of $(P_C)$ can theoretically be found by solving the Bellman equation considering  our problem $(P_C)$ as a perfect state information problem where the new state is $I_k$. If (\ref{sys_dyn}) and (\ref{eq_obs}) are linear, and $g_k$ and $g_T$ are quadratic in both the state and the control, the optimal policy is linear and can be computed in closed form. However, in the non-linear case, as the dimension of $I_k$ grows with time, $(P_C)$ is very often  intractable.
\subsection{Link with state estimation} \label{estimation}
If one supposes that $(\mu_0,\dots,\mu_{T-1})$ are constant, then only remains the minimization over $(\pi_0,\dots,\pi_{T})$ which gives a sequence of stochastic optimization problems, denoted by ${(P_E^k)}_{k \in \llbracket 0,T \rrbracket}$ that correspond to state estimation problems. For k $\in \llbracket 0,T \rrbracket$, $(P_E^k)$ is defined by :
\begin{equation*}
\begin{array}{rrclcc}
\displaystyle (P_E^k):  \underset{\pi_k}{\text{min}} & \multicolumn{3}{l}{f_C\left(\mathbb{E}\left[ (X_k-\widehat{X}_k){(X_k-\widehat{X}_k)}^T\right]\right)} \\
\textrm{s.t.} & \widehat{X}_{k} & = & \pi_k(I_k). \\
\end{array}
\end{equation*}

 If one chooses $f_C(\cdot)=\mathrm{tr}(\cdot)$ then:
 \begin{align*}
 f_C\left(\mathbb{E}\left[ (X_k-\widehat{X}_k){(X_k-\widehat{X}_k)}^T\right]\right)=\mathbb{E}\left[ {\Vert X_k-\widehat{X}_k\Vert}_2^2\right],
 \end{align*}
 and $(P_E^k)$ becomes the optimal filtering problem described in \cite{anderson_optimal_1979}
whose solution is known to be the conditional expectation of $X_k$ with respect to $I_k$ denoted by $\mathbb{E}[X_{k}\vert I_k ]$. If the equations (\ref{sys_dyn}) and (\ref{eq_obs}) are linear with independent gaussian disturbances then  $\mathbb{E}[X_{k}\vert I_k ]$ can be computed  exactly thanks to the recursive equations of the Kalman filter. Otherwise, such  exact equations do not exist and the problem becomes very hard.
Contrary to the min-max problem described in \cite{copp_nonlinear_2014}, in our case, when we combine the problems $(P_C)$ and  ${(P_E^k)}_{k \in \llbracket 0,T \rrbracket}$ to get $(P_{CE})$  the variables  $(U_0,\dots,U_{T-1})$ and $(\widehat{X}_0,\dots,\widehat{X}_T)$ are interestingly interdependent. 
\section{Tractable approximations of stochastic optimal control problems}\label{tract_app}
The optimal policy of $(P_{CE})$ denoted  by $(\mu^*_{0},\dots, \mu^*_{T-1},\pi^*_{0},\dots, \pi^*_{T})$ cannot be approached directly by space discretization of the information space because of its high dimension. So an approximation by a sub-optimal policy is proposed in this paper. First, before describing our approximation, we recall briefly a classification of stochastic control policies introduced in \cite{bar-shalom_dual_1974} and gives some example among existing methods. Then, we determine in which class our sub-optimal policy must be if we want to preserve the most important feature of $(P_{CE})$ that is to say the coupling between $U_{k-1}$ and $(\widehat{X}_k,\dots,\widehat{X}_T)$. Secondly, we explain how our approximation of
$(\mu^*_{0},\dots, \mu^*_{T-1},\pi^*_{0},\dots, \pi^*_{T})$, denoted by $(\mu^F_{0},\dots, \mu^F_{T-1},\pi^F_{0},\dots, \pi^F_{T})$ is computed.
\subsection{Classification of existing policies}

In \cite{bar-shalom_dual_1974}, four classes of stochastic control policies for fixed-end time are defined  according to the quantity of information used and the level of anticipation of the future. These classes of policies are defined as follow:
\begin{itemize}
{\setlength{\baselineskip}{1.1\baselineskip}
\item \textit{Open Loop} (OL) policies. In this case the control, $U_k$ for $k \in \llbracket 0,T-1 \rrbracket$,  depends only on the initial information $I_0$, the knowledge of dynamics (\ref{sys_dyn}) and of ${(p_{\xi_i})}_{\forall i \in \llbracket 0,T-1 \rrbracket}$. The sequence is determined once for all at time $k=0$ and never adapts itself to the available information. An application in robust path planning is described in \cite{blackmore_probabilistic_2010}
\item \textit{Feedback} (F) policies. In this class, $U_k$ depends on $I_k$, the dynamics (\ref{sys_dyn}), ${(p_{\xi_i})}_{\forall i \in \llbracket 0,T-1 \rrbracket}$, the observation equations (\ref{eq_obs}) up to time k and of ${(p_{\eta_i})}_{\forall i \in \llbracket 0,k \rrbracket}$. Unlike a OL-policy, a F-policy incorporates the current available information but never anticipates the fact that observations will be available at instants strictly greater than $k$. Many sub-optimal policies using \textit{Model Predictive Control} (MPC) combined with any estimator are F_policies because the fixed-time horizon optimization problems are solved with the initial condition being the current state estimate. MPC is used with particle filter in \cite{stahl_pf-mpc:_2011} and \cite{sehr_particle_2016}. Particle filter was already used to approximate stochastic control problems in  \cite{arnaud_doucet_sequential_2001}. This type of policies are reviewed in \cite{mayne_model_2014}. In \cite{copp_nonlinear_2014}, a policy combining worst case non-linear MPC and  a \textit{Moving Horizon Estimator} (MHE) into one global min-max problem is discussed. Still, this policy does not explicitly include knowledge of future observations so it remains a F-policy.
\item \textit{m-measurement feedback} (m-MF). In this class, $U_k$ depends on $I_k$, the dynamics (\ref{sys_dyn}), ${(p_{\xi_i})}_{\forall i \in \llbracket 0,T-1 \rrbracket}$, the observation equations (\ref{eq_obs}) up to time $k+m$ and of ${(p_{\eta_i})}_{\forall i \in \llbracket 0,k+m \rrbracket}$ with $m\leq T-k+1$. Similarly to  F-controls, m-MF- controls can adapt themselves to the current situation and also anticipate new observations up to $m$ instants after $k$. For example, Scenario-Based MPC \cite{mayne_model_2014} or Adaptive MPC \cite{mayne_model_2014} produces  m-MF policies. These controllers are said to be dual because, besides guiding the system to its initial goal, they also force the system to gain information about itself through state or parameter estimation. Examples of scenario based MPC are given in \cite{subramanian_economic_2015} and \cite{hanssen_scenario_2015}. Another example of a dual controller using particle filter and policy iterations is discussed in \cite{bayard_implicit_2008}.
\item \textit{Closed Loop} (CL) policies.  In this class, $U_k$ depends on $I_k$, the dynamics (\ref{sys_dyn}), ${(p_{\xi_i})}_{\forall i \in \llbracket 0,T-1 \rrbracket}$, the observation equations (\ref{eq_obs}) up to time T and of ${(p_{\eta_i})}_{\forall i \in \llbracket 0,T \rrbracket}$. This class is the extension of the m-MF class up to the final time T. Optimal policies obtained from Dynamic Programming belong to this class because  each  policy obtained from the backward Bellman equation minimizes a instantaneous cost plus a cost-to-go including all the future possible observations. We also suppose that $(\mu^*_{0},\dots, \mu^*_{T-1},\pi^*_{0},\dots, \pi^*_{T})$  belong to this class even if $f_C$ is not linear.
\par}
\end{itemize}
Considering this classification, our sub-optimal policy must belong at least to the m-MF class and ideally to the CL class. Indeed, the goal of our method is to get a control at time $k-1$, $U_{k-1}$, that reduces the estimation error made by $(\widehat{X}_k,\dots,\widehat{X}_T)$. We also know from equations (\ref{inf_vec}) and (\ref{func_constraint}) that  the estimator $(\widehat{X}_k,\dots,\widehat{X}_T)$ depends on the variables $(Z_0,\dots,Z_T)$. Besides, equations (\ref{sys_dyn}) and  (\ref{eq_obs}) show that the control $U_{k-1}$ cannot modify $(Z_0,\dots,Z_{k-1})$ but only $Z_k$ and by recursion the next observations $(Z_{k+1},\dots,Z_{T})$. Thus, our design of the control must incorporate the evolution of future observations thanks to equation (\ref{sys_dyn}) and (\ref{eq_obs}). If we consider $(Z_{k+1},\dots,Z_{T})$ then our sub-optimal policy is an CL policy. If we only include the evolution of $(Z_{k+1},\dots,Z_{k+m})$ for $m\leq T-k+1$, for computational reasons, then our policy is a m-MF one. In this paper, we described a version of our method that belongs to the CL class.

\subsection{Proposed particle approximation of the problem mixing control and state estimation}
Our approximation, $(\mu^F_{0},\dots, \mu^F_{T-1},\pi^F_{0},\dots, \pi^F_{T})$, is computed thanks to two separated ideas. First, we replace the term $f_C\left(C_k\right)$ by a term depending only on $X_k$ and $U_k$ removing the minimization over $\widehat{X}_k$ . Then, we approach the new problem by a sequence of deterministic problems solved online with a  technique similar to the one presented in \cite{sehr_particle_2016}. 
\subsubsection{Fisher approximation}
As said previously, $(\mu^F_{0},\dots, \mu^F_{T-1},\pi^F_{0},\dots, \pi^F_{T})$  must keep the coupling effect between the control and the state estimation. We recall from section \ref{sec_pce}  that the terms in the generalized cost that produce this effect are the terms $f_C\left(C_k\right)$. These terms also introduce a minimization over $(\widehat{X}_{0},\dots, \widehat{X}_T)$ without any other constraint that being a function of $I_k$, making $(P_{CE})$ impossible to approximate directly by several deterministic problems. The coupling disappears if one does this with a MPC-like technique and without modification in the cost. Indeed, if one transforms  $(P_{CE})$ in a deterministic problem fixing, for example, the disturbances to their mean, or with a Monte Carlo approximation then one does not look for policies anymore but for values so the constraints (\ref{func_constraint}) disappear. Then, $(\widehat{X}_{0}\cdots \widehat{X}_T)$ are unconstrained so, with, for instance, $f_C(\cdot)=\mathrm{tr}(\cdot)$, one finds $\forall k \in \llbracket 0,T \rrbracket$,  $X_k= \widehat{X}_k$. The computed value of $\widehat{X}_k$  is useless and the interesting terms also disappear. 
 To avoid this, we replace $C_k$ by ${\left(J_k\right)}^{-1}$ where $J_k$ is the \textit{Fisher Information Matrix (FIM)} which only depends on the current and previous states and  on the previous controls. Consequently, We have created a new stochastic optimization problem without optimization over $(\widehat{X}_{0},\dots, \widehat{X}_T)$. The new integral costs denoted by, $\tilde{g}_k^F$, and final cost denoted by, $\tilde{g}_T^F$, are then defined as follow $\forall k \in \llbracket 0,T-1 \rrbracket$:
 \begin{align}
 \tilde{g}^F_k(X_k,V_k,\xi_k) &={g}_k(X_k,U_k,\xi_k)+f_C\left({\left(J_k\right)}^{-1}\right),\\
\tilde{g}^F_T(X_T,V_T) &={g}_N(X_T)+f_C\left({\left(J_T\right)}^{-1}\right),
 \end{align}
where ${(J_k)}_{k \in \llbracket 0,T \rrbracket}$ is the FIM computed recursively  as in \cite{tichavsky_posterior_1998}.
 The new stochastic optimal control problem to solve is 
\begin{equation*}
\begin{array}{rrclcc}
\displaystyle (P_{CF}):\underset{\underset{ }{\mu_{0}\dots \mu_{T-1} }}{\text{min}} & \multicolumn{3}{l}{\mathbb{E}\left[\sum_{k=0}^{T-1} \tilde{g}^F_k(X_k,U_k,\xi_k) + \tilde{g}^F_T(X_T)\right]} \\
\textrm{s.t.} & X_{k+1} & = & f_k(X_{k},U_{k},{\xi}_{k}), \\
& Z_k & = & h_k(X_k,\eta_k), \\
& U_k & = & \mu_k(I_k), \; \forall k \in \llbracket 0,T-1 \rrbracket. \\
\end{array}
\end{equation*}
As the estimators are not included in the optimization problem anymore, we suppose that some estimators are computed outside of $({P}_{CF})$. Now, we have to justify why the coupling between the control and the state estimation still exists even if the estimators are not computed inside the optimization problem anymore.
We know from \cite{tichavsky_posterior_1998} that $J_k$ is invertible and for all non-biased estimator $ \widehat{X}_k$ of $X_k$, we have:
\begin{align*}
\mathbb{E}\left[(X_k-\widehat{X}_k){(X_k-\widehat{X}_k)}^T\right] \geq  J_{k}^{-1},
\end{align*}
where $\geq$ corresponds to a positive semi-definite inequality.
Moreover, let us assume that we choose an unbiased estimator $\widehat{X}_k$ whose covariance matrix $C_k$ tends to the inverse of the FIM when $k \rightarrow \infty$. Then if $f_C$ is continuous, minimizing $f_C\left({(J_k)}^{-1}\right)$ is close to minimizing $f_C(C_k)$ after a certain time . Thus the optimal policy of $(P_{CF})$ gives a control that almost minimizes   $f_C(C_k)$. In other words, the error made by  the estimator $\widehat{X}_k$ (in the sense of $f_C$) when estimating the optimal trajectory of  $({P}_{CF})$ is closed to be minimum. Consequently, the coupling between $U_{k-1}$ and $\widehat{X}_k$ still exists even if $\widehat{X}_k$ is removed from the optimization problem. This is also true for all the future estimators then one recovers the full coupling between $U_{k-1}$ and $(\widehat{X}_k,\dots,\widehat{X}_T)$.

\subsubsection{Particle approximation}
 The second idea consists in approximating $(P_{CF})$ by a Monte Carlo method. We use a set of particles and weights coming from a Particle Filter. Therefore, we suppose that, for  $ l \in \llbracket 0,T-1 \rrbracket$, the conditional density of $X_l$ w.r.t. $I_l$  denoted by $p(X_l\vert I_l)$ is represented by a set of N particles ${\left(\tilde{x}^{(i)}_{l}\right)}_{i\in\llbracket 1,N \rrbracket }$ and  weights ${\left(\omega^{(i)}_{l}\right)}_{i\in\llbracket 1,N \rrbracket }$.  
 This approximation is based on the fact that $p(X_l\vert I_l)$  is a sufficient statistics for classic problems with imperfect state information (\cite{bertsekas_dynamic_2005}, \cite{bayard_implicit_2008}) meaning that the policies can be considered as functions of $p(X_l\vert I_l)$ instead of $I_l$. Moreover, for computational reasons, we only use the $N_s<N$ most likely particles from the set ${\left(\tilde{x}^{(i)}_{l}\right)}_{i\in\llbracket 1,N \rrbracket }$.
We note that, in $(\tilde{P}_{CF}^l)$, the FIM is approximated with a Monte Carlo method. 

Our particle approximation of $(P_{CF})$ consists in solving a sequence of deterministic problems ${(\tilde{P}_{CF}^l)}_{l \in \llbracket 0,T-1 \rrbracket}$ defined by: 
\begin{equation*}
\begin{array}{rrclcc}
\displaystyle \underset{\underset{{x^{(i)}_{l+1}\cdots x^{(i)}_T }}{u_{l}\cdots u_{T-1}}}{\text{min}}\hspace{-0.3cm} & \multicolumn{3}{l}{\sum_{k=l}^{T-1}\sum_{i=1}^{N_s}\omega_l^{(i)} \left(\tilde{g}^F_k\left(x^{(i)}_{k},u_k,\xi_k \right) + \tilde{g}^F_T\left(x^{(i)}_{T}\right)\right)} \\
\textrm{s.t.} &\forall k& \in \llbracket l,T-1 \rrbracket,\\
&\forall i& \in \llbracket 1,N_s \rrbracket,&\\
& x^{(i)}_{l} & = & \tilde{x}^{(i)}_{l}, \\
& x^{(i)}_{k+1} & = &f_k(x^{(i)}_{k},u_{k},{\xi}^{(i)}_{k}).\\
\end{array}
\end{equation*}
Finally, for l  $\in \llbracket 0,T-1 \rrbracket$, we define our policy by:
\begin{align}
\mu^J_{l}\left(p(X_l\vert I_l)\right)&=u_l^*\label{control_app},\\
\pi^J_{l}(p(X_l\vert I_l))&=\sum_{i=1}^N \nolimits\omega_l^{(i)}\tilde{x}_l^{(i)} \label{estimation_app},\\
\pi^J_{T}(p(X_T\vert I_T))&=\sum_{i=1}^N \nolimits\omega_T^{(i)}\tilde{x}_T^{(i)}\label{estimation_app_fin}.
\end{align}
Equality (\ref{control_app}) means that we only apply the first optimal control found by solving $(\tilde{P}_{CF}^l)$. Equality (\ref{estimation_app}) and (\ref{estimation_app_fin}) mean that our estimator is $\mathbb{E}[X_{k}\vert I_k ]$ computed with a Monte Carlo method.	
Our feedback algorithm is summed up in Algorithm \ref{algo_fisher}.

\subsection{Application and Results}
\subsubsection{Description of our application}
\begin{figure}[h]
\begin{center}
\includegraphics[scale=0.4]{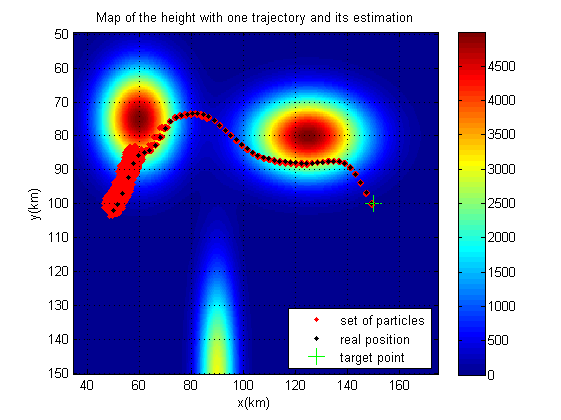}  
\includegraphics[scale=0.4]{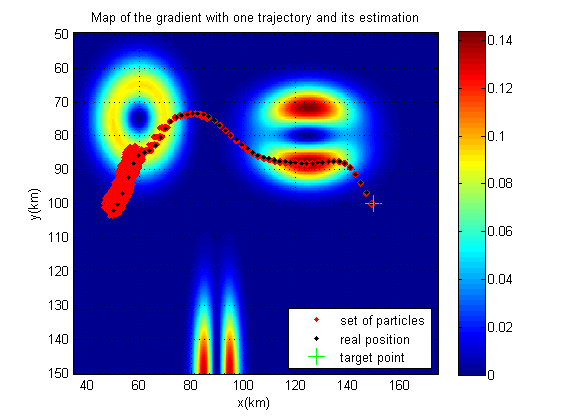}  
\end{center}
\caption{Plot of one trajectory obtained by fisher particle control and of the particles from the particle filter}
\label{fisher_traj}
\end{figure}
We applied this method to the guidance and localization of a drone by terrain-aided navigation. Our objective is to guide a drone in 3D from the unknown initial condition $X_0$ to a target point $x_{ta}$. To do so, we only measure the difference between the altitude of the drone and the altitude of the corresponding vertical point on the ground. More formally, at time $k$, the state $X_k$ is of dimension 6 and denoted $X_k={(x^1_k,x^2_k,x^3_k,v^1_k,v^2_k,v^3_k)}$ where $(x^1_k,x^2_k,x^3_k)$ stands for a 3D position and $(v^1_k,v^2_k,v^3_k)$ for a 3D speed. We suppose that (\ref{sys_dyn}) is linear i.e. $\forall k \in \llbracket 0,T-1 \rrbracket$:
\begin{align}
X_{k+1}&=FX_k+BU_k+\xi_k,\label{sys_lin},
\end{align}
Where $F$ and $B$ represent the discrete-time dynamic of a double integrator with a fixed time step dt.

To represent the observations made by the system we introduce $h_{map}: \mathbb{R}\times\mathbb{R} \longrightarrow \mathbb{R}$ which maps a horizontal position $(x^1,x^2)$ to the corresponding height on a terrain map. We suppose that $h_{map}$ is known but as it is often constructed from empirical data coming from a real terrain, it is highly non linear. Then the observation equation (\ref{eq_obs}) can be rewritten, $\forall k \in \llbracket 0,T \rrbracket$:
\begin{align}
Z_k=x^3_k-h_{map}(x^1_k,x^2_k)+ \eta_k.\label{eq_obs_map}
\end{align}
The challenge of this problem is to reconstruct a 6 dimensional state $X_k$ and, in particular, the horizontal position of the drone $(x^1_k,x^2_k)$ with a 1 dimensional observation. The main issue of this problem is that (\ref{sys_lin}) and (\ref{eq_obs_map}) may not be observable depending on the area the drone is flying over. Indeed, let us assume that the drone flies over a flat area then one measurement of height on the map correspond to a whole horizontal area  so the state estimation cannot be accurate. However, if the drone flies over a rough terrain, then one measurement of height matches a much smaller horizontal area and the state estimation can be more accurate. Therefore, the quantity that must be maximized is the gradient of $h_{map}$. Actually, from \cite{tichavsky_posterior_1998}, one can see that a quadratic term of this gradient appears in $J_k$ contains useful information to maintain the coupling between control and state estimation, as predicted in the previous part. 
\begin{algorithm}[h]
{\setlength{\baselineskip}{1\baselineskip}
\caption{Fisher Feedback Control}\label{algo_fisher}
\begin{algorithmic}[1]
\State Create a sample  of $N$ particles $\tilde{x}^{(i)}_{0}$ according to the law $\mathcal{N}(m_0,P_0)$ and initialize the weights $\omega_0^{(i)}$
\For{$l=0,\cdots,T-1$} 
\State Solve $(\tilde{P}_{CF}^l)$ starting from the set $\tilde{x}^{(i)}_{l}$ and the weights $\omega_k^{(i)}$.
\State Get a sequence of optimal control $u_l^*,\cdots,u_{T-1}^*$.
\State Draw  realizations of $\xi_l$, denoted by ${\xi}^{(i)}_{l}$.
\State Compute the \textit{a priori} set at time  $l$, ${\left(\tilde{x}^{(p,i)}_{l}\right)}_{i\in\llbracket 1,N \rrbracket }$, applying the dynamics (\ref{sys_dyn}) with control $u_l^*$ i.e: $  \tilde{x}^{(p,i)}_{l}=f_k(\tilde{x}^{(i)}_{l},u_{l}^*,{\xi}^{(i)}_{l})$.
\State Get the new observation $y_{l+1}$.
\State Compute the new weights ${\left(\omega^{(i)}_{l+1}\right)}_{i\in\llbracket 1,N \rrbracket }$.
\State  Compute the \textit{a posteriori} set ${\left(x^{(i)}_{l+1}\right)}_{i\in\llbracket 1,N \rrbracket }$  by re-sampling  the \textit{a priori} set ${\left(\tilde{x}^{(p,i)}_{l}\right)}_{i\in\llbracket 1,N \rrbracket }$.
\EndFor
\end{algorithmic}
\par}
\end{algorithm}
The desired online goal of our method in this application is then to estimate the state of the system and simultaneously design controls  that force the drone to fly over rough terrain so that the future estimation error diminishes. We also want the system to be guided precisely to the target $x_{ta}$, eventually.
Without state estimation improvement, we would like the drone to go in straight line to the target so we define the standard integral and final costs, $\forall k \in \llbracket 0,T-1 \rrbracket$, as follow:
\begin{align*}
&{g}_k(X_k,U_k,\xi_k)&=\alpha{\Vert U_k \Vert}^2_2,&\; &{g}_T(X_T)&= \gamma{\Vert X_T-x_{ta} \Vert}^2_2,&
\end{align*}
where $\alpha>0, \gamma >0$.
To generate the estimation improvement, we choose the coupling cost as follow:
\begin{align}
f_C\left({(J_k)}^{-1}\right)=\frac{\beta}{\mathrm{tr}(J_k)}, \forall k \in \llbracket 0,T \rrbracket, \label{cout_fisher}
\end{align}
where $\beta>0$. 
In section \ref{estimation}, we recalled that the natural cost would be $f_C\left({(J_k)}^{-1}\right)=\mathrm{tr}\left({(J_k)}^{-1}\right)$. However, in order to avoid matrix inverses in the resolution of $(\tilde{P}_{CJ}^l)$, we rather chose the cost defined in (\ref{cout_fisher}) that has the same monotony as the natural one in the matrix sense.
The parameters $(\alpha,\beta, \gamma) $ allow one to modify the behaviour of the system. If one wants to go faster to the target one can increase $\alpha$, on the contrary if one can afford to lose time and wants a more precise estimation then one can increase  $\beta$.
We have only applied our method on an artificial analytical but the final desired application is to use our method on real maps.
 \begin{figure}[h]
 \begin{center}
 \includegraphics[scale=0.45]{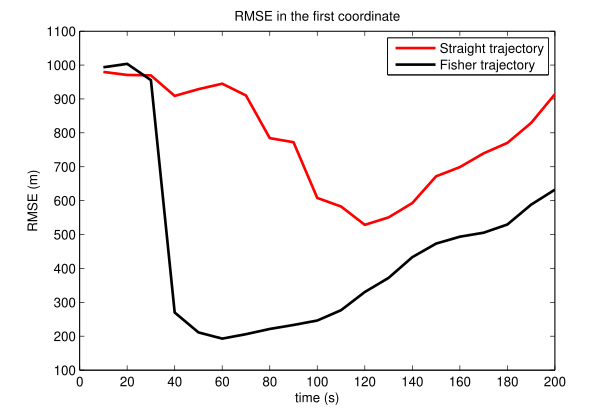}  
 \includegraphics[scale=0.45]{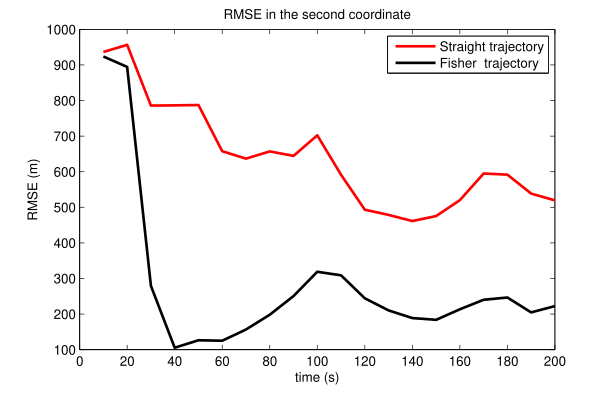} 
\end{center}
 \caption{Plot of the RMSE of a straight and a fisher trajectory in the $x^1$ and $x^2$ coordinates for 50 runs of algorithm \ref{algo_fisher} with $T=20$, $\mathrm{dt}=10s$ $N_s=100$ and $N=10000$}
 \label{RMSE}
\end{figure}
\subsubsection{Results} Figure \ref{fisher_traj} represents a simulated trajectory (black) in 2D of our drone   computed with the controls found by Algorithm \ref{algo_fisher} for one realization of the initial condition and of the disturbances. The figure also shows the particles (red) used to estimate the state of the  trajectory.  One can see that the set of particles tightens around the black trajectory. Other simulations have shown that it is not the case with a straight trajectory. Figure \ref{RMSE} compares the \textit{Root Mean Square Error} (RMSE) in $x^1$ and $x^2$  in the case of straight trajectories ($\beta=0$)  to the case of curved trajectories (large $\beta$) that creates coupling, for 50 runs of our algorithm. One can see that making a detour over the hills reduces highly the error made on the horizontal position of the drone compared to a standard trajectory, designed to go as fast as possible to the target. One can remark that in our example of map the RMSE in $x^1$ increases in both cases at the end of the runs. This is due to the ambiguity of our artificial terrain. One can also remark that our method allows the drone to avoid flat areas but not areas that would be non-flat and periodic. 

\section*{Conclusion}
This paper considers a stochastic optimal control problem combining state estimation and standard control designed to create dual effect. As this problem is intractable, a new approximation of the optimal control policy based on the FIM and a Particle Filter is proposed. Numerical results are given and show the efficiency of the whole method compared to the one without dual effect. In future works, from a theoretical point of view, we would like to evaluate the error made by solving $(P_{CF})$ with a fixed estimator instead of $(P_{CE})$. From an application point of view, we would like to apply the method on real maps and implement our method in a receding horizon way and a better Particle filter to decrease the number of particle needed and speed up the computations.

\bibliographystyle{plain}
\bibliography{articles_these}

\end{document}